\documentclass[12pt]{article}
\usepackage{amssymb,amsmath}
\usepackage{cases}
\usepackage{amsthm}
\usepackage{amsfonts}
\usepackage{graphicx}
\usepackage{cite,color,xcolor}
\usepackage[left=2.4cm,right=2.4cm,top=2.4cm,bottom=2.4cm]{geometry}
\usepackage[colorlinks,citecolor=blue,urlcolor=blue]{hyperref}
\usepackage[utf8]{inputenc}

\newtheorem{theorem}{Theorem}[section]
\newtheorem{corollary}{Corollary}[section]
\newtheorem{lemma}{Lemma}[section]
\newtheorem{proposition}{Proposition}[section]

\newtheorem{definition}{Definition}[section]
\newtheorem{remark}{Remark}[section]

\usepackage{txfonts}

\newcommand{\bal}{\begin{align}}
\newcommand{\bbal}{\begin{align*}}
\newcommand{\beq}{\begin{equation}}
\newcommand{\eeq}{\end{equation}}
\newcommand{\bca}{\begin{cases}}
\newcommand{\eca}{\end{cases}}

\newcommand{\pa}{\partial}
\newcommand{\fr}{\frac}

\newcommand{\De}{\Delta}

\newcommand{\ep}{\varepsilon}
\newcommand{\dd}{\mathrm{d}}

\newcommand{\R}{\mathbb{R}}
\newcommand{\T}{\mathbb{T}}

\newcommand{\les}{\lesssim}

\newcommand{\bi}{\Big}
\newcommand{\g}{\big}
\begin{document}
\bibliographystyle{plain}
\title{Non-uniform dependence on initial data for the Camassa--Holm equation in Besov spaces: Revisited}

\author{Jinlu Li$^{1}$, Yanghai Yu$^{2,}$\footnote{E-mail: lijinlu@gnnu.edu.cn; yuyanghai214@sina.com(Corresponding author); mathzwp2010@163.com} and Weipeng Zhu$^{3}$\\
\small $^1$ School of Mathematics and Computer Sciences, Gannan Normal University, Ganzhou 341000, China\\
\small $^2$ School of Mathematics and Statistics, Anhui Normal University, Wuhu 241002, China\\
\small $^3$ School of Mathematics and Big Data, Foshan University, Foshan, Guangdong 528000, China}

\date{\today}

\maketitle\noindent{\hrulefill}

{\bf Abstract:} In the paper, we revisit the uniform continuity properties of the data-to-solution map of the Camassa--Holm equation on the real-line case.
We show that the data-to-solution map of the Camassa--Holm equation is not uniformly continuous on the initial data in Besov spaces $B_{p, r}^s(\mathbb{R})$ with $s>\fr12$ and $1\leq p, r< \infty$, which improves the previous works \cite[Asian J. Math., 11 (2007)]{HMP07}, \cite[J. Differ. Equ., 269 (2020)]{Lyz} and \cite[J. Math. Fluid Mech., 23 (2021)]{Lwyz}. Furthermore, we present a strengthening of our previous work in \cite[J. Differ. Equ., 269 (2020)]{Lyz} and prove that the data-to-solution map for the Camassa--Holm equation is nowhere uniformly continuous in $B^s_{p,r}(\R)$ with $s>\max\{1+1/{p},3/2\}$ and $(p,r)\in [1,\infty]\times[1,\infty)$. The method applies also to the b-family of equations which contain the Camassa--Holm and Degasperis--Procesi equations.


\vskip0mm\noindent{\hrulefill}

\section{Introduction}

In this paper, we are concerned with the Cauchy problem for the classical Camassa--Holm (CH) equation
\begin{equation}\label{0}
\begin{cases}
u_t-u_{xxt}+3uu_x=2u_xu_{xx}+uu_{xxx}, \; &(x,t)\in \R\times\R^+,\\
u(x,t=0)=u_0(x),\; &x\in \R.
\end{cases}
\end{equation}
Here the scalar function $u = u(t, x)$ stands for the fluid velocity at time $t\geq0$ in the $x$ direction.

Setting $\Lambda^{-2}=(1-\pa^2_{xx})^{-1}$, then $\Lambda^{-2}f=G*f$ where $G(x)=\fr12e^{-|x|}$ is the kernel of the operator $\Lambda^{-2}$. Thus, we can transform the CH equation \eqref{0} equivalently into the following transport type equation
\begin{equation}\label{CH}
\begin{cases}
\partial_tu+u\pa_xu=\mathbf{P}(u), \; &(x,t)\in \R\times\R^+,\\
u(x,t=0)=u_0(x),\; &x\in \R,
\end{cases}
\end{equation}
where
\begin{equation}\label{CH1}
\mathbf{P}(u)=P(D)\left(u^2+\fr12(\pa_xu)^2\right)\quad\text{with}\quad P(D)=-\pa_x\Lambda^{-2}.
\end{equation}
The CH equation \eqref{0} was firstly proposed in the context of hereditary symmetries studied in \cite{Fokas} and then was derived explicitly as a water wave equation by Camassa--Holm \cite{Camassa}. Many aspects of the mathematical beauty of the CH equation have been exposed over the last two decades. Particularly, (CH) is completely integrable \cite{Camassa,Constantin-P} with a bi-Hamiltonian structure \cite{Constantin-E,Fokas} and infinitely many conservation laws \cite{Camassa,Fokas}. Also, it admits exact peaked
soliton solutions (peakons) of the form $u(x,t)=ce^{-|x-ct|}\;(c>0)$, which are orbitally stable \cite{Constantin.Strauss}. Another remarkable feature of the CH equation is the wave breaking phenomena: the solution remains bounded while its slope becomes unbounded in finite time \cite{Constantin,Escher2,Escher3}. It is worth mentioning that the peaked solitons present the characteristic for the travelling water waves of greatest height and largest amplitude and arise as solutions to the free-boundary problem for incompressible Euler equations over a flat bed, see Refs. \cite{Constantin-I,Escher4,Escher5,Toland} for the details.

Due to these interesting and remarkable features, the CH equation has attracted much attention as a class of integrable shallow water wave equations in recent twenty years. Its systematic mathematical study was initiated in a series of papers by Constantin and Escher, see \cite{Escher1,Escher2,Escher3,Escher4,Escher5}. After the CH equation was derived physically in the context of water waves, there are
a large amount of literatures devoted to studying the well-posedness of the Cauchy problem \eqref{0} (see
Molinet's survey \cite{Molinet}). Li and Olver \cite{li2000} proved that the Cauchy problem \eqref{0} is
locally well-posed with the initial data $u_0\in H^s(\R)$ with $s > 3/2$ (see also \cite{GB}). Danchin \cite{d1,d3} proved the local existence and uniqueness of strong solutions to \eqref{0} with initial data in $B^s_{p,r}$ if $(p,r)\in[1,\infty]\times[1,\infty), s>\max\g\{1+1/p, 3/2\g\}$ and $B^{3/2}_{2,1}$.  Meanwhile, he \cite{d1} only obtained the continuity of the solution map of \eqref{0} with respect to the initial data in the space $\mathcal{C}([0, T ];B^{s'}_{p,r})$ with any $s'<s$. Li-Yin \cite{Li-Yin1} proved the continuity of the solution map of \eqref{0} with respect to the initial data in the space $\mathcal{C}([0, T];B^{s}_{p,r})$ with $r<\infty$. In particular, they \cite{Li-Yin1} proved that the solution map of \eqref{0} is weak continuous with respect to initial data $u_0\in B^s_{p,\infty}$. For the endpoints, Danchin \cite{d3} obtained that the data-to-solution map is not continuous by using peakon solution, which implies the ill-posedness of \eqref{0} in $B^{3/2}_{2,\infty}$.  Guo-Liu-Molinet-Yin \cite{Guo-Yin} showed the ill-posedness of \eqref{0} in $B_{p,r}^{1+1/p}(\mathbb{R}\;\text{or}\; \mathbb{T})$ with $(p,r)\in[1,\infty]\times(1,\infty]$ (especially in $H^{3/2}$) by proving the norm inflation. Very recently, Guo-Ye-Yin \cite{Guo} obtained the ill-posedness for the CH equation in $B^{1}_{\infty,1}(\R)$ by proving the norm inflation. In our recent papers \cite{Li22,Li22-2}, we established the ill-posedness for \eqref{0} in $B^s_{p,\infty}(\mathbb{R})$ by proving the solution map to the CH equation starting from $u_0$ is discontinuous at $t = 0$ in the metric of $B^s_{p,\infty}(\mathbb{R})$.

From the PDE's point of view, it is crucial to know if an equation which models a physical phenomenon is well-posed in the
Hadamard's sense: existence, uniqueness, and continuous dependence of the solutions with respect to the initial data. In particular, the continuity of solution map is an important part of the well-posedness theory since the lack of continuous dependence would cause incorrect solutions or non meaningful solutions. Furthermore, the non-uniform continuity of data-to-solution map suggests that the local well-posedness cannot be established by the contraction mappings
principle since this would imply Lipschitz continuity for the solution map. After the phenomenon of non-uniform continuity for some dispersive equations was studied by Kenig et al. \cite{Kenig2001}, the issue of non-uniform dependence on the initial data has been a fascinating object of research in the recent past. Naturally, we may wonder which regularity assumptions are relevant for the initial data $u_0$ such that the Cauchy problem to \eqref{0} is not uniform dependent on initial data, namely, the dependence of solution on the initial data associated with this equation is not uniformly continuous. Himonas-Misio{\l}ek \cite{H-M} obtained the first result on the non-uniform dependence for \eqref{0} in $H^s(\T)$ with $s\geq2$ using explicitly constructed travelling wave solutions, which was sharpened to $s>\fr32$ by Himonas-Kenig \cite{H-K} on the real-line and Himonas-Kenig-Misio{\l}ek \cite{H-K-M} on the circle. We should mention that, non-uniform
continuity of the CH solution map in $H^1(\R\;\text{or}\;\T)$ was established by Himonas-Misio{\l}ek-Ponce \cite{HMP07}
by using traveling wave solutions. In our recent papers \cite{Lyz,Lwyz}, we proved the non-uniform dependence on initial data for \eqref{0} under both the framework of Besov spaces $B^s_{p,r}$ for $s>\max\big\{1+1/p, 3/2\big\}$ with $(p,r)\in[1,\infty]\times[1,\infty)$ and $B^{3/2}_{2,1}$. In particular, we established

\begin{theorem}[\cite{Lyz}]\label{th0}
Denote $U_R\equiv\g\{u_0\in B_{p,r}^s: \|u_0\|_{B^{s}_{p,r}}\leq R\g\}$ for any $R>0$. Assume that $(s,p,r)$ satisfies
\begin{align*}
s>\max\left\{\frac32,1+\frac1p\right\}\quad   \text{and}    \quad (p,r)\in [1,\infty]\times[1,\infty).
\end{align*}
Then the data-to-solution map of the Cauchy problem \eqref{CH}--\eqref{CH1}
\begin{equation*}
\mathbf{S}_t:\begin{cases}
U_R \rightarrow \mathcal{C}([0, T] ; B_{p, r}^{s}) \cap \mathcal{C}^1([0, T] ; B_{p, r}^{s-1}),\\
u_0\mapsto \mathbf{S}_t(u_0),
\end{cases}
\end{equation*}
is not uniformly continuous from any bounded subset $U_R$ in $B^s_{p,r}$ into $\mathcal{C}([0,T];B^s_{p,r})$. More precisely, there exists two sequences of solutions $\mathbf{S}_t(f_n+g_n)$ and $\mathbf{S}_t(f_n)$ such that
\bbal
&\|f_n\|_{B^s_{p,r}}\lesssim 1 \quad\text{and}\quad \lim_{n\rightarrow \infty}\|g_n\|_{B^s_{p,r}}= 0
\end{align*}
but
\bbal
\liminf_{n\rightarrow \infty}\|\mathbf{S}_t(f_n+g_n)-\mathbf{S}_t(f_n)\|_{B^s_{p,r}}\gtrsim t,  \quad \forall \;t\in[0,T_0],
\end{align*}
with small time $T_0$.
\end{theorem}
Ye-Yin-Guo \cite{Ye} proved the uniqueness and continuous dependence of the Camassa--Holm type equations in critical Besov spaces $B^{1+1/p}_{p,1}$ with $p\in[1,\infty)$. However, they do not consider the issue of non-uniform dependence on the initial data of the solution map. This is the starting point of this paper. We shall plan to investigate that whether or not the data-to-solution map of the Cauchy problem for the CH equation in the low regularity Besov spaces $B_{p,r}^{s}(\R)$ with $s>\fr12$ is uniformly continuous in this paper.

Our first result states as follows, which in particular improves the corresponding results in \cite{Lyz,Lwyz}.
\begin{theorem}\label{th-camassa}
Denote $U_R\equiv\g\{u_0\in B_{p,r}^s: \|u_0\|_{B^{s}_{p,r}}\leq R\g\}$ for any $R>0$. Assume that $(s,p,r)$ satisfies
\begin{align}\label{con}
(s,p,r)\in \left(\fr12,\infty\right)\times[1,\infty)\times [1,\infty).
\end{align}

\begin{description}
\item[(1)]\; If $s\in(1,\infty)$, the data-to-solution map of the Cauchy problem \eqref{CH}--\eqref{CH1}
\begin{equation*}
\mathbf{S}_t:\begin{cases}
U_R \rightarrow \mathcal{C}([0, T] ; B_{p, r}^{s}) \cap \mathcal{C}^1([0, T] ; B_{p, r}^{s-1}),\\
u_0\mapsto \mathbf{S}_t(u_0),
\end{cases}
\end{equation*}
is not uniformly continuous from any bounded subset $U_R$ in $B^s_{p,r}$ into $\mathcal{C}([0,T];B^s_{p,r})$.
More precisely, there exists two sequences of solutions $\mathbf{S}_t(f_n+g_n)$ and $\mathbf{S}_t(f_n)$ such that
\bbal
&\|f_n\|_{B^s_{p,r}}\lesssim 1 \quad\text{and}\quad \lim_{n\rightarrow \infty}\|g_n\|_{B^s_{p,r}}= 0
\end{align*}
but
\bbal
\liminf_{n\rightarrow \infty}\|\mathbf{S}_t(f_n+g_n)-\mathbf{S}_t(f_n)\|_{B^s_{p,r}}\gtrsim t,  \quad \forall \;t\in[0,T_0],
\end{align*}
with small time $T_0$.
  \item[(2)]\; If $s\in(\fr12,1]$, when the data-to-solution map $u_0\mapsto \mathbf{S}_t(u_0)$ of the Cauchy problem \eqref{CH}--\eqref{CH1} exists on $B^s_{p,r}(\R)$, (1) remains true.
\end{description}
\end{theorem}
\begin{remark}\label{rey}
Compared with \cite{Lyz,Lwyz}, Theorem \ref{th-camassa} is new since the regularity index is enlarged to $s>\fr12$. In this sense, Theorem \ref{th-camassa} improves the previous results in \cite{Lyz,Lwyz}.
\end{remark}
\begin{remark}\label{re4}
It should be mentioned that, $B^1_{2,2}(\R)=H^1(\R)$, thus Theorem \ref{th-camassa} covers the non-periodic result in \cite{HMP07}.
\end{remark}
\begin{remark}\label{re1}
The method we used in proving the Theorem \ref{th-camassa} is general and can be applied equally well to other related system, such as the Novikov equation
\begin{equation}\label{5}
\begin{cases}
u_t+u^2u_x=-\frac12(1-\pa^2_x)^{-1}u_x^3-\pa_x(1-\pa^2_x)^{-1}\left(\frac32uu^2_x+u^3\right),\\
u(x,t=0)=u_0(x).
\end{cases}
\end{equation}
By modifying the construction of initial data and repeating the procedure of the proof of Theorem \ref{th-camassa},
we can obtain the similar results for the Novikov equation. Since the process is standard, we leave the details  to the interested readers.
\end{remark}
\begin{remark}\label{re10}
We should emphasize that, for the low regularity case, it is not possible to proceed as in \cite{Lwyz} due to the initial constructions: specifically for more general critical case $s=1+1/p$, we encounter the failure of the uniform bound owing to the rougher initial data. By the subtle modification of the strategies used in \cite{Lwyz}, we decompose the solution maps as
\bbal
&\mathbf{S}_{t}(\underbrace{f_n+g_n}_{=~u^n_0})=\underbrace{\mathbf{S}_{t}(u^n_0)-u^n_0+tu^n_{0}\pa_xu^n_{0}}_{=~\mathbf{I}_1(u^n_0)}+f_n+g_n-tu^n_{0}\pa_xu^n_{0}
\end{align*}
and
\bbal
\mathbf{S}_{t}(f_n)=\underbrace{\mathbf{S}_{t}(f_n)-f_n+tf_n\pa_x f_n}_{=~\mathbf{I}_2(f_n)}+f_n-tf_n\pa_xf_n.
\end{align*}
We expect that $u^n_{0}\pa_xu^n_{0}$ brings us the term $g_n\partial_xf_n$ which still plays an essential role since it would not small when $n$ is large enough. To deal with the errors estimates $\mathbf{I}_2(f_n)$ and $\mathbf{I}_1(u^n_{0})$ in $B^{s}_{p,r}$, the common point is that, we need to establish the $L^\infty$-estimation of the solution map and finetuned $C^{0,1}$-estimation of $\mathbf{S}_{t}(u_0)-u_0$ since the embedding $B^{s}_{p,r}(\R)\hookrightarrow L^\infty(\R)$ is no longer available. The novelty of the proof of our result lies in the new observation, $B^{s}_{p,r}$ and $L^\infty$ possess different level of regularity, which allows us to appropriately define these new initial data sequences $f_n$ and $g_n$. Based on the suitable choice of $f_n$ and $g_n$, we shall prove that for a short time $t\in(0,1]$ and for some $\varepsilon_n\rightarrow0$ as $n\rightarrow\infty$
\begin{itemize}
  \item $\liminf\limits_{n\rightarrow \infty}\|g_n\pa_xf_n\|_{B^s_{p,r}}\gtrsim 1,$
  \item $\|\mathbf{I}_2(f_n)\|_{B^{s}_{p,r}}+\|\mathbf{I}_1(u^n_{0})\|_{B^{s}_{p,r}}\les t^2+\varepsilon_n,$
\end{itemize}
which means that the data-to-solution map is not uniformly continuous.
\end{remark}
Our second main result of this paper is a stronger version of Theorem \ref{th0}.
\begin{theorem}[{\bf Nowhere uniformly continuous}]\label{th4}
Assume that $(s,p,r)$ satisfies
\begin{align}\label{con-camassa}
s>\max\left\{\frac32,1+\frac1p\right\}\quad   \text{and}    \quad (p,r)\in [1,\infty]\times[1,\infty).
\end{align}
{\bf Then the data-to-solution map of the Cauchy problem \eqref{CH}--\eqref{CH1} is nowhere uniformly continuous from $B^s_{p,r}$ into $\mathcal{C}([0,T];B^s_{p,r})$} in the following sense: For any $u_0$ in $B^s_{p,r}(\R)$, there exists some $B^s_{p,r}$-neighbourhood $\mathcal{V}_{u_0}$ of $u_0$ such that for any $v_0\in \mathcal{V}_{u_0}$, the system \eqref{CH}--\eqref{CH1} has a unique solution $V\in \mathcal{C}([0, T] ; B_{p, r}^{s}) \cap \mathcal{C}^1([0, T] ; B_{p, r}^{s-1})$. Furthermore, the data-to-solution map
\begin{equation*}
\mathbf{S}_t:\begin{cases}
\mathcal{V}_{u_0} \rightarrow \mathcal{C}([0, T] ; B_{p, r}^{s}) \cap \mathcal{C}^1([0, T] ; B_{p, r}^{s-1}),\\
v_0\mapsto V=\mathbf{S}_t(v_0),
\end{cases}
\end{equation*}
is not uniformly continuous on $\mathcal{V}_{u_0}$.
More precisely, for any $u_0\in B^s_{p,r}(\R)$, there exists two initial sequences $v^1_{0,n}=u_0+\mathbf{f}_n$ and $v^2_{0,n}=u_0+\mathbf{f}_n+\mathbf{g}_n$ such that two sequences of solutions $\mathbf{S}_t(v^1_{0,n})$ and $\mathbf{S}_t(v^2_{0,n})$ satisfying
\bbal
&\|\mathbf{f}_n\|_{B^s_{p,r}}\lesssim 1 \quad\text{and}\quad \lim_{n\rightarrow \infty}\|v^2_{0,n}-v^1_{0,n}\|_{B^s_{p,r}}=\lim_{n\rightarrow \infty}\|\mathbf{g}_n\|_{B^s_{p,r}}= 0
\end{align*}
but
\bbal
\liminf_{n\rightarrow \infty}\|\mathbf{S}_t(v^2_{0,n})-\mathbf{S}_t(v^1_{0,n})\|_{B^s_{p,r}}\gtrsim t,  \quad \forall \;t\in[0,T_0],
\end{align*}
with small time $T_0$.
\end{theorem}
\begin{remark}\label{relyz}
Theorem \ref{th0} implies that the data-to-solution map $\mathbf{S}_t$ has not the property to be uniformly continuous on bounded sets while Theorem \ref{th4} means that for any $u_0\in B^s_{p,r}$ the restriction $\mathbf{S}_t|_{\mathcal{V}_{u_0}}$ is not uniformly continuous. In this sense, Theorem \ref{th4} improves the previous results in \cite{Lyz}.
\end{remark}
\begin{remark}
We also mention that, using a geometric approach, Inci obtained a series of nowhere-uniform continuity results includes the b-family of equations \cite{Inci16}, the two component b-family of equations \cite{Inci22}, the hyperelastic rod equation \cite{Inci19}, and the incompressible Euler equation \cite{Incili}. Bourgain and Li \cite{Bou} showed that the data-to-solution map for the incompressible Euler equations is nowhere-uniform continuity in $H^s(\R^d)$ with $s\geq0$ by using an idea of localized Galilean boost. It is worth mentioning that although our result is obtained by pure analytic method which is more or less elementary, it seems that the method is more simpler and very robust, and can be applied equally well to other related systems as mentioned above.
\end{remark}
\begin{remark}
Holm and Staley \cite{03A} introduced the following b-family equation (see \cite{03B,DP,Li22dp} etc.):
\begin{align}\label{b-family}
\begin{cases}
\pa_tu+u\pa_xu=-\pa_x(1-\pa^2_x)^{-1}\left(\frac{b}{2}u^2
+\frac{3-b}{2}u^2_x\right), &\quad (t,x)\in \R^+\times\R,\\
u(0,x)=u_0(x), &\quad x\in \R.
\end{cases}
\end{align}
\end{remark}
It should be emphasised that the Camassa--Holm equation corresponds to $b = 2$ and Degasperis--Procesi equation corresponds to $b = 3$.
Since the concrete values of the parameter $b$ have no impact on the proof of Theorem \ref{th4}, thus Theorem \ref{th4} also holds for the b-family equations \eqref{b-family}. Due to the fact that the Besov space $B^s_{2,2}$ coincides with the Sobolev space
$H^s$, Theorem \ref{th4} covers the previous result given by Inci \cite{Inci16} who proved that the corresponding solution map of b-family of equations in the Sobolev spaces $H^s(\R)$ with $s>\fr32$ is nowhere locally uniformly continuous.

When $b=3$, \eqref{b-family} reduces to the Degasperis--Procesi (DP) equation
\begin{equation}\label{DP}
\begin{cases}
\partial_tu+u\pa_xu=-\frac32\pa_x(1-\pa^2_x)^{-1}u^2, \; &(x,t)\in \R\times\R^+,\\
u(x,t=0)=u_0(x),\; &x\in \R.
\end{cases}
\end{equation}
In particular, we have
\begin{theorem}\label{th5}
Assume that $(s,p,r)$ satisfies
\begin{align*}
s>1+\frac1p, \ (p,r)\in [1,\infty]\times [1,\infty) \quad \mbox{or} \quad
s=\frac{1}{p}+1, \ (p,r)\in [1,\infty]\times\{1\}.
\end{align*} Then the DP equation \eqref{DP} is nowhere uniformly continuous from $B^s_{p,r}$ into $\mathcal{C}([0,T];B^s_{p,r})$.
\end{theorem}
\begin{remark}\label{re2}
Following the procedure in the proof of Theorem \ref{th4} with suitable modification, we can prove Theorem \ref{th5}. Here we will omit the details.
\end{remark}

\quad\noindent\textbf{Organization of our paper.} In Section \ref{sec2}, we list some notations and known results and recall some Lemmas which will be used in the sequel. In Section \ref{sec3}, we establish some technical Lemmas and Propositions and then prove Theorem \ref{th-camassa}. In Section \ref{sec5}, we establish two crucial Propositions and then prove Theorem \ref{th4}.

\section{Preliminaries}\label{sec2}
\subsection{Notation}\label{subsec21}
We will use the following notations throughout this paper.
Given a Banach space $X$, we denote its norm by $\|\cdot\|_{X}$. For $I\subset\R$, we denote by $\mathcal{C}(I;X)$ the set of continuous functions on $I$ with values in $X$. Sometimes we will denote $L^p(0,T;X)$ by $L_T^pX$. We will also define the Lipschitz space $C^{0,1}$ using the norm $\|f\|_{C^{0,1}}=\|f\|_{L^\infty}+\|\pa_xf\|_{L^\infty}$.
The symbol $a\lesssim b$ means that there is a uniform positive constant $C$ independent of $a$ and $b$ such that $a\leq Cb$.
Let us recall that for all $u\in \mathcal{S}'$, the Fourier transform $\mathcal{F}u$, also denoted by $\widehat{u}$, is defined by
$
\mathcal{F}u(\xi)=\widehat{u}(\xi)=\int_{\R}e^{-\mathrm{i}x\xi}u(x)\dd x$ for any $\xi\in\R.
$
\subsection{Littlewood-Paley analysis}\label{subsec22}
Next, we will recall some facts about the Littlewood-Paley decomposition, the nonhomogeneous Besov spaces and their some useful properties (see \cite{B} for more details).

There exists a couple of smooth functions $(\chi,\varphi)$ valued in $[0,1]$, such that $\chi$ is supported in the ball $\mathcal{B}= \{\xi\in\mathbb{R}:|\xi|\leq \frac 4 3\}$, and $\varphi$ is supported in the ring $\mathcal{C}= \{\xi\in\mathbb{R}:\frac 3 4\leq|\xi|\leq \frac 8 3\}$. Moreover, for any $\xi\in \R$, there holds
\begin{align*}
&\chi(\xi)+\sum_{j\geq0}\varphi(2^{-j}\xi)=1 \quad \mbox{ and }\quad
\fr13\leq\chi^2(\xi)+\sum_{j\geq0}\varphi^2(2^{-j}\xi)\leq1.
\end{align*}
We should emphasize that the fact $\varphi(\xi)\equiv 1$ for $4/3\leq |\xi|\leq 3/2$ will be used in the sequel.

For every $u\in \mathcal{S'}(\mathbb{R})$, the inhomogeneous dyadic blocks ${\Delta}_j$ are defined as follows
\begin{numcases}{\Delta_ju=}
0, & if $j\leq-2$;\nonumber\\
\chi(D)u=\mathcal{F}^{-1}(\chi \mathcal{F}u), & if $j=-1$;\nonumber\\
\varphi(2^{-j}D)u=\mathcal{F}^{-1}\g(\varphi(2^{-j}\cdot)\mathcal{F}u\g), & if $j\geq0$.\nonumber
\end{numcases}
The inhomogeneous low-frequency cut-off operator $S_{j}$ is defined by
$$
S_j u=\sum_{-1\leq q\leq j-1}{\Delta}_qu.
$$
With this we give the definition of nonhomogeneous Besov space $B^{s}_{p,r}(\R)$.
\begin{definition}[see \cite{B}]
Let $s\in\mathbb{R}$ and $(p,r)\in[1, \infty]^2$. The nonhomogeneous Besov space $B^{s}_{p,r}(\R)$ is defined by
\begin{align*}
B^{s}_{p,r}(\R):=\Big\{f\in \mathcal{S}'(\R):\;\|f\|_{B^{s}_{p,r}(\mathbb{R})}<\infty\Big\},
\end{align*}
where
\begin{numcases}{\|f\|_{B^{s}_{p,r}(\mathbb{R})}=}
\left(\sum_{j\geq-1}2^{sjr}\|\Delta_jf\|^r_{L^p(\mathbb{R})}\right)^{1/r}, &if $1\leq r<\infty$,\nonumber\\
\sup_{j\geq-1}2^{sj}\|\Delta_jf\|_{L^p(\mathbb{R})}, &if $r=\infty$.\nonumber
\end{numcases}
\end{definition}
\begin{remark}
 It is easy to verify that the nonhomogeneous Besov space $B^{s}_{2,2}(\mathbb{R})$ coincides with the classical nonhomogeneous Sobolev space $H^{s}(\mathbb{R})$. It should be emphasized that the following embedding will be often used implicity:
$$B^s_{p,q}(\R)\hookrightarrow B^t_{p,r}(\R)\quad\text{for}\;s>t\quad\text{or}\quad s=t,1\leq q\leq r\leq\infty.$$
 \end{remark}

\begin{lemma}[see \cite{B}]\label{densely}
If $r$ is finite, then for any $u \in B^s_{p,r}(\R)$, we have
\bbal
&\lim_{j\to\infty}\|S_ju-u\|_{B^{s}_{p,r}(\R)}=0.
\end{align*}
\end{lemma}
We recall the following facts involving the product estimate and interpolation inequality, which will be used in the later of the paper.
\begin{lemma}[see \cite{B}]\label{le1}
Let $(p,r)\in[1, \infty]^2$ and $s>0$. For any $u,v \in B^s_{p,r}(\R)\cap L^\infty(\R)$, we have
\bbal
&\|uv\|_{B^{s}_{p,r}(\R)}\leq C(\|u\|_{B^{s}_{p,r}(\R)}\|v\|_{L^\infty(\R)}+\|v\|_{B^{s}_{p,r}(\R)}\|u\|_{L^\infty(\R)}).
\end{align*}
In particular, for $s>\frac1p$, due to the fact $B^{s}_{p,r}(\R)\hookrightarrow L^\infty(\R)$, then we have
\begin{align*}
&\|uv\|_{B^{s}_{p,r}(\R)}\leq C\|u\|_{B^{s}_{p,r}(\R)}\|v\|_{B^{s}_{p,r}(\R)}.
\end{align*}
\end{lemma}
\begin{lemma}[see \cite{B}]\label{le2} If $s_{1}$ and $s_{2}$ are real numbers such that $s_{1}<s_{2},$
$\theta \in(0,1)$, and $(p, r)\in[1, \infty]^2$, then we have
\begin{align*}
\|u\|_{B_{p, 1}^{\theta s_{1}+(1-\theta) s_{2}}} \leq \frac{C}{s_{2}-s_{1}}\left(\frac{1}{\theta}+\frac{1}{1-\theta}\right)\|u\|_{B_{p, \infty}^{s_{1}}}^{\theta}\|u\|_{B_{p, \infty}^{s_{2}}}^{1-\theta}.
\end{align*}
\end{lemma}

\begin{lemma}[see \cite{B}]\label{le3}
Let $(p,r)\in[1, \infty]^2$ and $\sigma\geq-\min\g\{\frac1p, 1-\frac1p\g\}$. Assume that $f_0\in B^\sigma_{p,r}(\R)$, $g\in L^1([0,T]; B^\sigma_{p,r}(\R))$ and
\begin{numcases}{\pa_{x}\mathbf{u}\in}
L^1([0,T]; B^{\sigma-1}_{p,r}(\R)), & if $\sigma>1+\fr1p$ \mbox{or}\ $\sigma=1+\fr1p,\; r=1$;\nonumber\\
L^1([0,T]; B^{\sigma}_{p,r}(\R)), & if $\sigma=1+\fr1p,\; r>1$;\nonumber\\
L^1([0,T]; B^{1/p}_{p,\infty}(\R)\cap L^\infty(\R)), & if $\sigma<1+\fr1p$.\nonumber
\end{numcases}
If $f\in L^\infty([0,T]; B^\sigma_{p,r}(\R))\cap \mathcal{C}([0,T]; \mathcal{S}'(\R))$ solves the following linear transport equation:
\begin{equation*}
\quad \partial_t f+\mathbf{u}\pa_xf=g,\quad \; f|_{t=0} =f_0.
\end{equation*}
\begin{enumerate}
\item Then there exists a constant $C=C(p,r,\sigma)$ such that the following statement holds
\begin{equation*}
\|f(t)\|_{B^\sigma_{p,r}(\R)}\leq e^{CV(t)} \Big(\|f_0\|_{B^\sigma_{p,r}(\R)}+\int_0^t e^{-CV(\tau)} \|g(\tau)\|_{B^\sigma_{p,r}(\R)}\mathrm{d}\tau\Big),
\end{equation*}
where
\begin{numcases}{V(t)=}
\int_0^t \|\pa_x\mathbf{u}(\tau)\|_{B^{\sigma-1}_{p,r}(\R)}\mathrm{d}\tau,\; &if $\sigma>1+\fr1p$ \mbox{or} $\sigma=1+\fr1p,\; r=1$;\nonumber\\
\int_0^t \|\pa_x\mathbf{u}(\tau)\|_{B^{\sigma}_{p,r}(\R)}\mathrm{d}\tau, \; &if $\sigma=1+\fr1p,\; r>1$;\nonumber\\
\int_0^t \|\pa_x\mathbf{u}(\tau)\|_{B^{1/p}_{p,\infty}(\R)\cap L^\infty(\R)}\mathrm{d}\tau,\; &if $\sigma<1+\fr1p$.\nonumber
\end{numcases}
\item If $\sigma>0$, then there exists a constant $C=C(p,r,\sigma)$ such that the following statement holds
\begin{align*}
\|f(t)\|_{B^\sigma_{p,r}}&\leq \|f_0\|_{B^\sigma_{p,r}}+C\int_0^t\|g(\tau)\|_{B^\sigma_{p,r}}\mathrm{d}\tau\\
&\quad+C\int^t_0\Big(\|f(\tau)\|_{B^\sigma_{p,r}}\|\pa_x\mathbf{u}(\tau)\|_{L^\infty}+\|\pa_x\mathbf{u}(\tau)\|_{B^{\sigma-1}_{p,r}}\|\pa_x f(\tau)| |_{L^\infty}\Big)\mathrm{d}\tau.
\end{align*}
\end{enumerate}
\end{lemma}
\begin{lemma}[see Lemma 3.26 in \cite{B}]\label{l22}
Let $(p,r)\in[1, \infty]^2$, $s>1$ and $u_0\in B^s_{p,r}(\R)$.
Assume that $u\in L^\infty([0,T]; B^s_{p,r}(\R))$ solves \eqref{CH}--\eqref{CH1}.
Then there exists a constant $C=C(s,p)$ and a universal constant $C'$ such that for all $t\in[0,T]$, we have
\begin{align}
&\|u(t)\|_{B^s_{p,r}(\R)}\leq \|u_0\|_{B^s_{p,r}(\R)}\exp\left(C\int_0^t \|u(\tau)\|_{C^{0,1}(\R)}\mathrm{d}\tau\right),\label{c1}\\
&\|u(t)\|_{C^{0,1}(\R)}\leq \|u_0\|_{C^{0,1}(\R)}\exp\left(C\int_0^t \|\pa_xu(\tau)\|_{L^{\infty}(\R)}\mathrm{d}\tau\right).\label{c2}
\end{align}
\end{lemma}
\begin{corollary}\label{corr1}
Let $(p,r)\in[1, \infty]^2$, $s>1$ and $u_0\in B^s_{p,r}\cap C^{0,1}$.
Assume that $u\in L^\infty([0,T]; B^s_{p,r}(\R))$ solves \eqref{CH}--\eqref{CH1}.
There exists a constant $C>0$ such that for all $t\in(0,\min\{1,1/(2C\|u_0\|_{C^{0,1}}\})]$, we have
\begin{align*}
&\|u(t)\|_{C^{0,1}}\leq C\|u_0\|_{C^{0,1}},\\
&\|u(t)\|_{B^s_{p,r}}\leq C\|u_0\|_{B^s_{p,r}}.
\end{align*}
\end{corollary}
\begin{proof} Setting
\bbal
\lambda(t):= \|u_0\|_{C^{0,1}}\exp\left(C\int_0^t \sup_{s\in[0,\tau]}\|\pa_xu(s)\|_{L^{\infty}(\R)}\dd \tau\right) \quad\text{with}\quad \lambda(0):= \|u_0\|_{C^{0,1}},
\end{align*}
then from \eqref{c2}, one has
\bbal
\frac{\dd}{\dd t}\lambda(t)\leq C\lambda^2(t)\quad\Leftrightarrow\quad -\frac{\dd}{\dd t}\left(\frac{1}{\lambda(t)}\right)\leq C.
\end{align*}
Solving the above directly yields for $t\in(0,\min\{1,1/(2C\|u_0\|_{C^{0,1}}\})]$
\bbal
\sup_{\tau\in[0,t]}\|u(\tau)\|_{C^{0,1}}\leq \frac{\|u_0\|_{C^{0,1}}}{1-Ct\|u_0\|_{C^{0,1}}}\leq 2\|u_0\|_{C^{0,1}}.
\end{align*}
From \eqref{c1}, we also have
\begin{align*}
&\sup_{\tau\in[0,t]}\|u(\tau)\|_{B^s_{p,r}}\leq C\|u_0\|_{B^s_{p,r}}.
\end{align*}
This completes the proof of Corollary \ref{corr1}.\end{proof}

From Lemma \ref{l22}, we have the local well-posedness for the CH equation in Besov spaces.
\begin{lemma}[see \cite{d1,d3}]\label{well}
Assume that $(s,p,r)$ satisfies \eqref{con-camassa} and for any initial data $u_0\in B_{p,r}^s$. Then there exists some $T=T(\|u_0\|_{B_{p,r}^s})>0$ such that the CH equation \eqref{CH}--\eqref{CH1} has a unique solution $\mathbf{S}_{t}(u_0)\in \mathcal{C}([0,T];B^s_{p,r})$. Moreover, we have
\begin{align*}
\|\mathbf{S}_{t}(u_0)\|_{B^s_{p,r}}\leq C\|u_0\|_{B_{p,r}^s}.
\end{align*}
\end{lemma}

\section{Proof of Theorem \ref{th-camassa}}\label{sec3}
\subsection{Technical Lemmas}\label{subsec3}
Firstly, we need to introduce smooth, radial cut-off functions to localize the frequency region. Precisely,
let $\widehat{\phi}\in \mathcal{C}^\infty_0(\mathbb{R})$ be an even, real-valued and non-negative function on $\R$ and satisfy
\begin{numcases}{\widehat{\phi}(\xi)=}
1,&if $|\xi|\leq \frac{1}{4}$,\nonumber\\
0,&if $|\xi|\geq \frac{1}{2}$.\nonumber
\end{numcases}
Motivated by \cite{Lyz}, we establish the following crucial lemmas which will be used later on.
\begin{lemma}[see \cite{Lyz}]\label{l31} For any $p\in[1,\infty]$, then there exists a positive constant $M$ such that
\begin{align*}
\liminf_{\lambda\rightarrow + \infty}\left\|\phi^2(x)\cos \left(\lambda x\right)\right\|_{L^p}\geq M.
\end{align*}
\end{lemma}
From now on, to prove Theorem \ref{th-camassa}, we denote $\delta=(1-\eta)p$, where
\begin{numcases}{\eta=}
\fr12,&if $s\geq1$,\nonumber\\
\fr12\left(s-\fr12\right),&if $\frac{1}{2}<s<1$.\nonumber
\end{numcases}
\begin{lemma}\label{l32} Let $s\in\R$ and $p\in[1,\infty)$.
Define the high frequency function $f_n$ by
\bbal
&f_n=2^{-n(s+1-\eta)}\phi\left(2^{-\delta n}x\right)\sin \left(\frac{17}{12}2^nx\right),\quad n\gg1.
\end{align*}
Then for any $\sigma\in\R$, there exists a positive constant $C=C(\phi)$ such that
\bal
&\|f_n\|_{L^\infty}\leq C2^{-n(s+1-\eta)},\label{m1}\\
&\|\pa_xf_n\|_{L^\infty}\leq C2^{-n(s-\eta)},\label{m2}\\
&\|f_n\|_{B^\sigma_{p,r}}\leq C2^{n(\sigma-s)}.\label{m3}
\end{align}
\end{lemma}
\begin{proof} \eqref{m1} and \eqref{m2} are obvious. Next we shall prove \eqref{m3}. Due to the simple fact $\sin\theta=\frac{\mathrm{i}}{2} (e^{-\mathrm{i} \theta}-e^{\mathrm{i} \theta})$, we deduce easily that
\begin{align*}
\mathcal{F}\left[\phi\left(2^{-\delta n} x\right) \sin \left(\frac{17}{12} 2^{n} x\right)\right]
=\;&\int_{\R} e^{-\mathrm{i} x \xi} \phi\left(2^{-\delta n} x\right) \sin \left(\frac{17}{12} 2^{n} x\right) \dd x \\
=\;&\frac{\mathrm{i}}{2} \int_{\R}\left(e^{-\mathrm{i} x\left(\xi+\frac{17}{12} 2^{n}\right)}-e^{-\mathrm{i} x\left(\xi-\frac{17}{12} 2^{n}\right)}\right) \phi\left(2^{-\delta n} x\right) \dd x \\
=\;&\frac{\mathrm{i}}{2}\cdot2^{\delta n} \int_{\R}\left(e^{-\mathrm{i} 2^{\delta n} x\left(\xi+\frac{17}{12} 2^{n}\right)}-e^{-\mathrm{i} 2^{\delta n} x\left(\xi-\frac{17}{12} 2^{n}\right)}\right) \phi(x) \dd x\\
=\;&\frac{\mathrm{i}}{2}\cdot2^{\delta n} \left[\widehat{\phi}\left(2^{\delta n} \left(\xi+\frac{17}{12} 2^{n}\right)\right)-\widehat{\phi}\left(2^{\delta n} \left(\xi-\frac{17}{12} 2^{n}\right)\right)\right],
\end{align*}
which is contained in the union of two balls centered at $\pm\frac{17}{12} 2^{n}$
and radius $2^{-(1+\delta n)}$, and hence in a single dyadic annulus, namely,
\bal\label{zml}
\mathrm{supp} \ \widehat{f}_n&\subset \left\{\xi\in\R: \ \frac{17}{12}2^n-\fr1{2^{1+\delta n}}\leq |\xi|\leq \frac{17}{12}2^n+\fr1{2^{1+\delta n}}\right\}\\
&\subset \left\{\xi\in\R: \ \frac{33}{24}2^{n}\leq |\xi|\leq \frac{35}{24}2^{n}\right\}.\nonumber
\end{align}
Recalling that $\varphi\equiv 1$ for $\frac43\leq |\xi|\leq \frac32$, we have
\begin{numcases}{\Delta_j(f_n)=}
f_n, &if $j=n$,\nonumber\\
0, &if $j\neq n$.\nonumber
\end{numcases}
By the definitions of Besov space, we deduce that
\bbal
\|f_n\|_{B^\sigma_{p,r}}&=2^{n\sigma}\|f_n\|_{L^p}=2^{n\sigma}2^{-n(s+1-\eta)}\left\|\phi\left(2^{-(1-\eta)p n}x\right)\sin \left(\frac{17}{12}2^nx\right)\right\|_{L^p}\\
&\leq2^{n\sigma}2^{-n(s+1-\eta)}\left\|\phi\left(2^{-(1-\eta)p n}x\right)\right\|_{L^p}\\
&\leq 2^{n(\sigma-s)}\|\phi\|_{L^p}.
\end{align*}
Thus we finished the proof of Lemma \ref{l32}.
\end{proof}

\begin{lemma}\label{l33} Let $s\in\R$ and $p\in[1,\infty)$. Define the low frequency function $g_n$ by
\bbal
&g_n=2^{-n}\phi\left(2^{-\delta n}x\right),\quad n\gg1.
\end{align*}
Then for any $\sigma\in\R$,  there exist two positive constants $C=C(\phi)$ and $c=c(\phi)$ such that
\bal
&\|g_n\|_{L^\infty}\leq C2^{-n}, \label{m4}\\
&\|g_n\|_{B^\sigma_{p,r}}\leq C2^{-\eta n}, \label{m5}\\
&\liminf_{n\rightarrow \infty}\|g_n\pa_xf_n\|_{B^s_{p,\infty}}\geq c. \label{m6}
\end{align}
\end{lemma}
\begin{proof} \eqref{m4} is obvious. Easy computations give that
$$
\mathcal{F}\left[\phi\left(2^{-\delta n} x\right)\right]=2^{\delta n} \widehat{\phi}\left(2^{\delta n} \xi\right),
$$
which implies
\bal\label{zml2}
\mathrm{supp} \ \widehat{g}_n\subset \left\{\xi\in\R: \ 0\leq |\xi|\leq \fr1{2^{1+\delta n}}\right\},
\end{align}
then, we have
\bbal
\widehat{\Delta_jg_n}=\varphi(2^{-j}\xi)\widehat{g}_n(\xi)\equiv0\quad\text{for}\quad j\geq0,
\end{align*}
equivalently,
\bbal
{\Delta_jg_n}\equiv0\quad\text{for}\quad j\geq0.
\end{align*}
Noticing that the definitions of $g_n$ and the Besov space, one obtains that
\bbal
\|g_n\|_{B^\sigma_{p,r}}&=2^{-(n+\sigma)}\bi\|\De_{-1}\phi\left(2^{-\delta n}x\right)\bi\|_{L^p}
\leq C2^{-(n+\sigma)}\bi\|\phi\left(2^{-\delta n}x\right)\bi\|_{L^p}
\leq C2^{-\eta n}\|\phi\|_{L^p},
\end{align*}
which gives \eqref{m5}.

Notice that \eqref{zml} and \eqref{zml2}, then we have
\bbal
\mathrm{supp}\ \widehat{g_n\pa_xf_n}\subset \left\{\xi\in\R: \ \frac{17}{12}2^n-\fr1{2^{\delta n}}\leq |\xi|\leq \frac{17}{12}2^n+\fr1{2^{\delta n}}\right\},
\end{align*}
which implies
\begin{numcases}{\Delta_j\g(g_n\pa_xf_n\g)=}
g_n\pa_xf_n, &if $j=n$,\nonumber\\
0, &if $j\neq n$.\nonumber
\end{numcases}
Thus, we have
\bal\label{yh}
\|g_n\pa_xf_n\|_{B^s_{p,\infty}}&=2^{ns}\|\De_{n}\g(g_n\pa_xf_n\g)\|_{L^p}=2^{ns}\|g_n\pa_xf_n\|_{L^p}.
\end{align}
By the definitions of $f_n$ and $g_n$, we obtain
\bbal
2^{ns}\left(g_n\pa_xf_n\right)&=I_1+I_2,
\end{align*}
where
\bbal
&I_1:=\frac{17}{12}\cdot2^{-(1-\eta)n}\phi^2\left(2^{-\delta n}x\right)\cos \left(\frac{17}{12}2^nx\right),\\
&I_2:=2\cdot2^{-(2-\eta)n}\cdot2^{-\delta n}\phi\left(2^{-\delta n}x\right)\pa_x\phi\left(2^{-\delta n}x\right)\sin \left(\frac{17}{12}2^nx\right).
\end{align*}
By the change of variables $y=2^{-\delta n}x$, then
\bbal
&\|I_1\|_{L^p}\geq \left\|\phi^2\left(x\right)\cos \left(\frac{17}{12}2^{(1+\delta) n}x\right)\right\|_{L^p},\\
&\|I_2\|_{L^p}\leq 2\cdot2^{-(1+\delta) n}\left\|\phi\left(x\right)\pa_x\phi\left(x\right)\sin \left(\frac{17}{12}2^{(1+\delta) n}x\right)\right\|_{L^p}\leq C(\phi)2^{-(1+\delta) n}.
\end{align*}
Using the above estimates and the triangle inequality, we obtain from \eqref{yh} that
\bal\label{z}
\|g_n\pa_xf_n\|_{B^s_{p,\infty}}&\geq \left\|\phi^2\left(x\right)\cos \left(\frac{17}{12}2^{(1+\delta) n}x\right)\right\|_{L^p}-C(\phi)2^{-(1+\delta) n}.
\end{align}
Using Lemma \ref{l31} to \eqref{z} enables us finish the proof of Lemma \ref{l33}.
\end{proof}
\begin{lemma}\label{l41} The data-to-solution map $u_0\mapsto \mathbf{S}_t(u_0)$ of the Cauchy problem \eqref{CH}--\eqref{CH1} satisfies that for $t\in(0,1]$
\bbal
\|\mathbf{S}_{t}(u_0)-u_0\|_{L^\infty}\leq Ct\|u_0\|^2_{C^{0,1}}.
\end{align*}
\end{lemma}
\begin{proof}
By the fundamental theorem of calculus in the time variable and Corollary \ref{corr1}, we have
\bbal
\|\mathbf{S}_{t}(u_0)-u_0\|_{L^\infty}&\leq\int^t_0\|\pa_\tau \mathbf{S}_{\tau}(u_0)\|_{L^\infty} \dd\tau
\nonumber\\&\leq \int^t_0\|\mathbf{S}_\tau(u_0)\pa_x\mathbf{S}_\tau(u_0)\|_{L^\infty}\dd \tau+\int^t_0\|\mathbf{P}\left(\mathbf{S}_\tau(u_0)\right)\|_{L^\infty}\dd \tau\nonumber\\
&\les\int^t_0\|\mathbf{S}_\tau(u_0)\|_{L^\infty}\|\pa_x\mathbf{S}_\tau(u_0)\|_{L^\infty}\dd \tau
+ \int^t_0\bi\|\pa_xG\ast\bi((\mathbf{S}_\tau(u_0))^2+\fr12(\pa_x\mathbf{S}_{t}(u_0))^2\bi)\bi\|_{L^\infty}\dd \tau
\nonumber\\
&\les\int^t_0\|\mathbf{S}_\tau(u_0)\|^2_{C^{0,1}}\dd \tau\nonumber\\
&\les t\|u_0\|^2_{C^{0,1}}.
\end{align*}
This completes the proof of Lemma \ref{l41}.
\end{proof}
The following proposition claims that the initial data sequence $u_0=f_n$ can approximate to the solution map $\mathbf{S}_t(f_n)$.
\begin{proposition}\label{pro1}
Assume that $u_0\in C^{0,1}$ and $\|u_0\|_{B^s_{p,r}}\lesssim 1$ with $(s,p,r)$ satisfying \eqref{con}. Then the data-to-solution map $u_0\mapsto \mathbf{S}_t(u_0)$ of the Cauchy problem \eqref{CH}--\eqref{CH1} satisfies that for $t\in(0,1]$
\bal
&\|\mathbf{S}_{t}(u_0)-u_0\|_{B^{s}_{p,r}}\les t\bi(\|u_0\|_{B^{s+1}_{p,r}}(\|u_0\|^2_{C^{0,1}}+\|u_0\|_{L^\infty})+\|u_0\|_{C^{0,1}}\|u_0\|_{B^{s+\ep}_{p,r}}\bi),\label{w1}\\
&\|\mathbf{S}_{t}(u_0)-u_0\|_{B^{s+1}_{p,r}}\les t\bi(\|u_0\|_{B^{s+2}_{p,r}}(\|u_0\|^2_{C^{0,1}}+\|u_0\|_{L^\infty})+\|u_0\|_{C^{0,1}}\|u_0\|_{B^{s+1}_{p,r}}\bi),\label{w2}
\end{align}
here and in what follows, we take
\begin{numcases}{\ep=}
\fr12,&if $\frac{1}{2}<s<1$,\nonumber\\
\in\left(0,\fr12\right),&if $s=1$,\nonumber\\
0,&if $s>1$.\nonumber
\end{numcases}
\end{proposition}
\begin{proof} By Corollary \ref{corr1}, we know that the solution map $\mathbf{S}_{t}(u_0)\in \mathcal{C}([0,T];B^s_{p,r})$ and has common lifespan $T\thickapprox1$. Moreover, there holds for $k\in\{0,1,2\}$
\bbal
&\|\mathbf{S}_{t}(u_0)\|_{L^\infty_T(B^{s+k}_{p,r})}\leq C\|u_0\|_{B^{s+k}_{p,r}},\quad s>1,\\
&\|\mathbf{S}_{t}(u_0)\|_{L^\infty_T(B^{s}_{p,r})}\leq C\|u_0\|_{B^{s+\ep}_{p,r}},\quad \frac{1}{2}<s\leq1.
\end{align*}
By the fundamental theorem of calculus in the time variable, we obtain
\bbal
\|\mathbf{S}_{t}(u_0)-u_0\|_{B^{s}_{p,r}}&\les \int^t_0\|\pa_\tau \mathbf{S}_{\tau}(u_0)\|_{B^{s}_{p,r}}\dd \tau\\
&\les\int^t_0\|\mathbf{S}_\tau(u_0)\pa_x\mathbf{S}_\tau(u_0)\|_{B^{s}_{p,r}}\dd \tau+\int^t_0\|\mathbf{P}\left(\mathbf{S}_\tau(u_0)\right)\|_{B^{s}_{p,r}}\dd \tau\\
&\les\int^t_0\|\mathbf{S}_\tau(u_0)\|_{L^\infty}\|\mathbf{S}_\tau(u_0)\|_{B^{s+1}_{p,r}}\dd \tau+\int^t_0\|\left(\mathbf{S}_\tau(u_0)\right)^2+\left(\pa_x\mathbf{S}_\tau(u_0)\right)^2\|_{B^{s-1}_{p,r}}\dd \tau\\
&\les\|u_0\|_{B^{s+1}_{p,r}}\int^t_0(\|\mathbf{S}_\tau(u_0)-u_0\|_{L^\infty}+\|u_0\|_{L^\infty})\dd \tau+\|\mathbf{S}_\tau(u_0)\|_{B^{s+\ep}_{p,r}}\int^t_0\|\mathbf{S}_\tau(u_0)\|_{C^{0,1}}\dd \tau\\
&\les t\bi(\|u_0\|_{B^{s+1}_{p,r}}(\|u_0\|^2_{C^{0,1}}+\|u_0\|_{L^\infty})+\|u_0\|_{C^{0,1}}\|u_0\|_{B^{s+\ep}_{p,r}}\bi),
\end{align*}
where we have used
$$\|\mathbf{P}(\mathbf{S}_{\tau}(u_0))\|_{B^s_{p,r}}\les\left\|\left(\mathbf{S}_\tau(u_0)\right)^2
+\left(\pa_x\mathbf{S}_\tau(u_0)\right)^2\right\|_{B^{s-1}_{p,r}}\les\|u_0\|_{C^{0,1}}\|u_0\|_{B^{s+\ep}_{p,r}}$$
and the fact from Lemma \ref{l41}
\bbal
\|\mathbf{S}_{t}(u_0)-u_0\|_{L^\infty}\les \|u_0\|^2_{C^{0,1}}.
\end{align*}
Similarly, we have
\bbal
\|\mathbf{S}_{t}(u_0)-u_0\|_{B^{s+1}_{p,r}}&\les \int^t_0\|\pa_\tau \mathbf{S}_{\tau}(u_0)\|_{B^{s+1}_{p,r}}\dd \tau\\
&\les\int^t_0\|\mathbf{S}_\tau(u_0)\pa_x\mathbf{S}_\tau(u_0)\|_{B^{s+1}_{p,r}}\dd \tau+\int^t_0\|\mathbf{P}\left(\mathbf{S}_\tau(u_0)\right)\|_{B^{s+1}_{p,r}}\dd \tau\\
&\les\int^t_0\|\mathbf{S}_\tau(u_0)\|_{L^\infty}\|\mathbf{S}_\tau(u_0)\|_{B^{s+2}_{p,r}}\dd \tau+\int^t_0\|\left(\mathbf{S}_\tau(u_0)\right)^2+\left(\pa_x\mathbf{S}_\tau(u_0)\right)^2\|_{B^{s}_{p,r}}\dd \tau\\
&\les\|u_0\|_{B^{s+2}_{p,r}}\int^t_0(\|\mathbf{S}_\tau(u_0)-u_0\|_{L^\infty}+\|u_0\|_{L^\infty})\dd \tau+\int^t_0\|\mathbf{S}_{t}(u_0)\|_{C^{0,1}}\|\mathbf{S}_{t}(u_0)\|_{B^{s+1}_{p,r}}\dd \tau\\
&\les t\bi(\|u_0\|_{B^{s+2}_{p,r}}(\|u_0\|^2_{C^{0,1}}+\|u_0\|_{L^\infty})+\|u_0\|_{C^{0,1}}\|u_0\|_{B^{s+1}_{p,r}}\bi).
\end{align*}
Then we complete the proof of Proposition \ref{pro1}.
\end{proof}
To obtain the non-uniformly continuous dependence property for the CH equation, we need to prove that the initial data sequence $u_0=f_n+g_n$ can not approximate to the solution map $\mathbf{S}_t(f_n+g_n)$.
\begin{proposition}\label{pro2}
Assume that $u_0\in C^{0,1}$ and $\|u_0\|_{B^s_{p,r}}\lesssim 1$ with $(s,p,r)$ satisfying \eqref{con}. Then the data-to-solution map $u_0\mapsto \mathbf{S}_t(u_0)$ of the Cauchy problem \eqref{CH}--\eqref{CH1} satisfies that for $t\in(0,1]$
\bal\label{l4}
\|\mathbf{S}_{t}(u_0)-u_0+tu_0\pa_x u_0\|_{B^{s}_{p,r}}&\leq Ct^2\mathbf{E}(u_0)+C\|u_0\|_{C^{0,1}}\|u_0\|_{B^{s+\ep}_{p,r}},
\end{align}
where we denote
\bbal
&\mathbf{E}(u_0)=\|u_0\|_{B^{s+1}_{p,r}}\|u_0\|^2_{C^{0,1}}
+\|u_0\|_{L^\infty}\|u_0\|_{B^{s+2}_{p,r}}\left(\|u_0\|^2_{C^{0,1}}+\|u_0\|_{L^\infty}
\right).
\end{align*}
\end{proposition}
\begin{proof} By the fundamental theorem of calculus in the time variable, from \eqref{CH}, we have
\begin{align*}
&\|\mathbf{S}_t(u_0)-u_0+tu_0\pa_x u_0\|_{B^s_{p,r}}\\
\leq&~ \int^t_0\|\partial_\tau \mathbf{S}_{\tau}(u_0)+u_0\pa_x u_0\|_{B^s_{p,r}} \dd \tau\\
\leq&~ \int^t_0\|\mathbf{S}_{\tau}(u_0)\partial_x\mathbf{S}_{\tau}(u_0)-u_0\partial_xu_0\|_{B^s_{p,r}}\dd \tau+\int^t_0\|\mathbf{P}(\mathbf{S}_{\tau}(u_0))\|_{B^s_{p,r}} \dd \tau
\\
\lesssim&~ \|u_0\|_{B^{s+1}_{p,r}}\int^t_0\|\mathbf{S}_{\tau}(u_0)-u_0\|_{L^\infty}\dd\tau
+\|u_0\|_{L^\infty}\int^t_0\|\mathbf{S}_{\tau}(u_0)-u_0\|_{B_{p, r}^{s+1}}\dd\tau+\|u_0\|_{C^{0,1}}\|u_0\|_{B^{s+\ep}_{p,r}}
\\
\lesssim&~t\left(\|u_0\|_{B^{s+1}_{p,r}}\|\mathbf{S}_{\tau}(u_0)-u_0\|_{L^\infty_tL^\infty}
 +\|\mathbf{S}_{\tau}(u_0)-u_0\|_{L^\infty_tB_{p, r}^{s+1}}\|u_0\|_{L^\infty}\right)+\|u_0\|_{C^{0,1}}\|u_0\|_{B^{s+\ep}_{p,r}},
\end{align*}
where we have used
\begin{align*}
 \|\mathbf{S}_{\tau}(u_0)\partial_x\mathbf{S}_{\tau}(u_0)-u_0\partial_xu_0\|_{B^s_{p,r}}
 \les&~\|(\mathbf{S}_{\tau}(u_0)-u_0)(\mathbf{S}_{\tau}(u_0)+u_0)\|_{B_{p, r}^{s+1}}\\
  \lesssim&~ \|\mathbf{S}_{\tau}(u_0)-u_0\|_{L^\infty}\|u_0\|_{B_{p, r}^{s+1}}+\|\mathbf{S}_{\tau}(u_0)-u_0\|_{B_{p, r}^{s+1}}\|\mathbf{S}_{\tau}(u_0)+u_0\|_{L^\infty}\\
  \lesssim&~ \|u_0\|_{B^{s+1}_{p,r}}\|\mathbf{S}_{\tau}(u_0)-u_0\|_{L^\infty}+\|\mathbf{S}_{\tau}(u_0)-u_0\|_{B_{p, r}^{s+1}}\|u_0\|_{L^\infty}.
   \end{align*}
Combining Lemma \ref{l41} and Proposition \ref{pro1} enable us to complete the proof of Proposition \ref{pro2}.
\end{proof}
With Propositions \ref{pro1}--\ref{pro2} in hand, we can prove Theorem \ref{th-camassa}.

{\bf Proof of Theorem \ref{th-camassa}.} We set $u^n_0=f_n+g_n$ and compare the solution $\mathbf{S}_{t}(u^n_0)$ and $\mathbf{S}_{t}(f_n)$.
Obviously, we have
\bbal
\|u^n_0-f_n\|_{B^s_{p,r}}=\|g_n\|_{B^s_{p,r}}\leq C2^{-\eta n}\quad \Rightarrow\quad
\lim_{n\to\infty}\|u^n_0-f_n\|_{B^s_{p,r}}=0.
\end{align*}
Notice that
\bbal
&\mathbf{S}_{t}(\underbrace{f_n+g_n}_{=~u^n_0})=\underbrace{\mathbf{S}_{t}(u^n_0)-u^n_0+tu^n_{0}\pa_xu^n_{0}}_{=~\mathbf{I}_1(u^n_0)}+f_n+g_n-tu^n_{0}\pa_xu^n_{0}\nonumber\\
&\mathbf{S}_{t}(f_n)=\underbrace{\mathbf{S}_{t}(f_n)-f_n+tf_n\pa_xf_n}_{=~\mathbf{I}_2(f_n)}+f_n-tf_n\pa_xf_n\quad\text{and}\nonumber\\
&u^n_{0}\pa_xu^n_{0}-f_n\pa_xf_n=g_n\pa_xf_n+u^n_{0}\pa_xg_n,
\end{align*}
using the triangle inequality and Proposition \ref{pro1}, we deduce that
\bal\label{h0}
\quad \ &\|\mathbf{S}_{t}(f_n+g_n)-\mathbf{S}_{t}(f_n)\|_{B^s_{p,r}}
=\|\mathbf{I}_1(u^n_0)-\mathbf{I}_2(f_n)+g_n-t\big(g_n\pa_xf_n+u^n_{0}\pa_xg_n\big)\|_{B^s_{p,r}}\nonumber\\
\geq&~t\|g_n\pa_xf_n\|_{B^s_{p,r}}-t\|u^n_{0}\pa_xg_n\|_{B^s_{p,r}}-\|\mathbf{I}_1(u^n_0)\|_{B^s_{p,r}}-\|\mathbf{I}_2(f_n)\|_{B^s_{p,r}}
-\|g_n\|_{B^s_{p,r}}\nonumber\\
\geq&~ t\|g_n\pa_xf_n\|_{B^s_{p,r}}-C2^{-\eta n}-C\|\mathbf{I}_1(u^n_0)\|_{B^s_{p,r}}-C\|\mathbf{I}_2(f_n)\|_{B^s_{p,r}},
\end{align}
where we have used
\bbal
&\|u^n_{0}\pa_xg_n\|_{B^s_{p,r}}\les \|f_n\|_{L^\infty}\|g_n\|_{B^{s+1}_{p,r}}+\|f_n\|_{B^s_{p,r}}\|\pa_xg_n\|_{L^\infty}+ \|g_n\|_{L^\infty}\|g_n\|_{B^{s+1}_{p,r}}\les2^{-n}.
\end{align*}
Using Proposition \ref{pro2} with $u_0=f_n$  and Lemma \ref{l32}, yields
\bal\label{h1}
&\|\mathbf{I}_2(f_n)\|_{B^{s}_{p,r}}\les t^2+2^{-n(s-\eta-\ep)}.
\end{align}
Using Proposition \ref{pro2} with $u_0=f_n+g_n$ and Lemmas \ref{l32}-\ref{l33}, yields
\bal\label{h2}
\|\mathbf{I}_1(u^n_0)\|_{B^s_{p,r}}\les t^2+2^{-n(s-\eta-\ep)}+2^{-(1-\ep)n}.
\end{align}
Using \eqref{h1} and \eqref{h2}, then \eqref{h0} reduces to
\bal\label{h3}
\liminf_{n\rightarrow \infty}\|\mathbf{S}_{t}(f_n+g_n)-\mathbf{S}_{t}(f_n)\|_{B^s_{p,r}}
\geq&~ t\liminf_{n\rightarrow \infty}\|g_n\pa_xf_n\|_{B^s_{p,r}}-Ct^{2}.
\end{align}
Hence, it follows from \eqref{h3} and Lemma \ref{l33} that
\bbal
\liminf_{n\rightarrow \infty}\|\mathbf{S}_t(f_n+g_n)-\mathbf{S}_t(f_n)\|_{B^s_{p,r}}\gtrsim t\quad\text{for} \ t \ \text{small enough}.
\end{align*}
This completes the proof of Theorem \ref{th-camassa}.{\hfill $\square$}

\section{Proof of Theorem \ref{th4}}\label{sec5}

\quad In this section, our aim is to proving Theorem \ref{th4}.
Firstly, we need to establish the useful Lemma.
\begin{lemma}\label{cor1}
Let $1 \leq p, r \leq \infty$ and $s>\max \{1+1 / p, 3 / 2\}$. Suppose that we are given
$$
(u, v) \in\left(L^{\infty}([0, T] ; B_{p, r}^{s}) \cap \mathcal{C}([0, T] ; B_{p, r}^{s-1})\right)^{2}
$$
two solutions of \eqref{CH}--\eqref{CH1}  with initial data $u_{0}, v_{0} \in B_{p, r}^{s}$. Letting $w:=u-v$, then we have, for every $t \in[0, T]$ and some constant $C=C(s, p, r)$
\begin{align}
&\|w(t)\|_{B_{p, r}^{s-1}} \leq\left\|w_{0}\right\|_{B_{p, r}^{s-1}} \exp \left(C \int_{0}^{t}\left\|(u,v)\left(\tau\right)\right\|_{B_{p, r}^{s}} \dd \tau\right),\label{cor1-1}\\
&\|w(t)\|_{B_{p, r}^{s}} \leq\left(\left\|w_{0}\right\|_{B_{p, r}^{s}}+C \int_{0}^{t}
\left\|w\right\|_{B_{p, r}^{s-1}}\left\|\partial_{x} v\right\|_{B_{p, r}^{s}} \dd\tau\right) \exp \left(C \int_{0}^{t}\left\|(u,v)\left(\tau\right)\right\|_{B_{p, r}^{s}} \dd \tau\right).\label{cor1-2}
\end{align}
\end{lemma}
\begin{proof} \eqref{cor1-1} is a straightforward result of Proposition 3.23 in \cite{B}. It remains to prove \eqref{cor1-2}. It is obvious that $w$ solves the transport equation
$$
\partial_{t} w+u \partial_{x} w=-w \partial_{x} v+\mathcal{B}(w, u+v)\quad\text{with}\quad
\mathcal{B}:(f, g) \mapsto P(D)\left(f g+\frac{1}{2} \partial_{x} f \partial_{x} g\right).
$$
According to Lemma \ref{le3} and Lemma \ref{le1}, the following inequality holds true:
\begin{align*}
\|w(t)\|_{B_{p, r}^{s}} &\leq\left\|w_{0}\right\|_{B_{p, r}^{s}} +C \int_{0}^{t}\left\|u\right\|_{B_{p, r}^{s}}\|w\|_{B_{p, r}^{s}} \dd\tau+C \int_{0}^{t}
\left(\left\|w \partial_{x} v\right\|_{B_{p, r}^{s}}+\|\mathcal{B}(w, u+v)\|_{B_{p, r}^{s}}\right) \dd\tau\\
&\leq\left\|w_{0}\right\|_{B_{p, r}^{s}} +C \int_{0}^{t}\left(\left\|u\right\|_{B_{p, r}^{s}}+\left\|v\right\|_{B_{p, r}^{s}}\right)\|w\|_{B_{p, r}^{s}} \dd\tau+C \int_{0}^{t}
\left\|w\right\|_{B_{p, r}^{s-1}}\left\|\partial_{x} v\right\|_{B_{p, r}^{s}} \dd\tau.
\end{align*}
Applying Gronwall's inequality yields \eqref{cor1-2} and thus completes the proof of Lemma \ref{cor1}.
\end{proof}

We define the high-low frequency functions $\mathbf{f}_n$ and $\mathbf{g}_n$ by
\bal\label{fg}
\mathbf{f}_n=2^{-ns}\phi(x)\sin \left(\frac{17}{12}2^nx\right)\quad\text{and}\quad \mathbf{g}_n=\frac{12}{17}2^{-n}\phi(x),\quad n\gg1.
\end{align}
Then we set $$v^{n,\omega}_{0,m_n}:=(\mathbf{f}_n+\omega \mathbf{g}_n)(x-m_n)\quad\text{with}\quad \omega\in\{0,1\},$$
where $m_n$ can be chosen later.

\begin{proposition}\label{pro5-1}
Under the assumptions of Theorem \ref{th4}, we have
\bbal
\|\mathbf{S}_{t}(u_0+v^{n,\omega}_{0,m_n})-\mathbf{S}_{t}(S_nu_0+v^{n,\omega}_{0,m_n})\|_{B^{s}_{p,r}}\leq C\|(\mathrm{Id}-S_n)u_0\|_{B^{s}_{p,r}}.
\end{align*}
\end{proposition}
\begin{proof} The local well-posedness result (see Lemma \ref{well}) tells us that $\mathbf{S}_{t}(u_0+v^{n,\omega}_{0,m_n}),\mathbf{S}_{t}(S_nu_0+v^{n,\omega}_{0,m_n})\in \mathcal{C}([0,T];B^s_{p,r})$ and has common lifespan $T\thickapprox1$. Moreover, there holds
\bal\label{s}
\|\mathbf{S}_{t}(u_0+v^{n,\omega}_{0,m_n}),\mathbf{S}_{t}(S_nu_0+v^{n,\omega}_{0,m_n})\|_{L^\infty_T(B^s_{p,r})}\leq C.
\end{align}
Using Lemma \ref{cor1} with $u=\mathbf{S}_{t}(u_0+v^{n,\omega}_{0,m_n})$ and $v=\mathbf{S}_{t}(S_nu_0+v^{n,\omega}_{0,m_n})$ yields that
\bal\label{yhy1}
\|\mathbf{S}_{t}(u_0&+v^{n,\omega}_{0,m_n})-\mathbf{S}_{t}(S_nu_0+v^{n,\omega}_{0,m_n})\|_{B^{s-1}_{p,r}}
\leq C\|(\mathrm{Id}-S_n)u_0\|_{B^{s-1}_{p,r}}
\leq C2^{-n}\|(\mathrm{Id}-S_n)u_0\|_{B^s_{p,r}}
\end{align}
and
\bbal
&\|\mathbf{S}_{t}(u_0+v^{n,\omega}_{0,m_n})-\mathbf{S}_{t}(S_nu_0+v^{n,\omega}_{0,m_n})\|_{B^{s}_{p,r}}\\
\leq&~ C\left(\|(\mathrm{Id}-S_n)u_0\|_{B^{s}_{p,r}}
+2^n\int^t_0\|\mathbf{S}_{\tau}(u_0+v^{n,\omega}_{0,m_n})-\mathbf{S}_{\tau}(S_nu_0+v^{n,\omega}_{0,m_n})\|_{B^{s-1}_{p,r}}\dd \tau\right)\quad \text{by}\; \eqref{yhy1}\\
\leq&~ C\|(\mathrm{Id}-S_n)u_0\|_{B^{s}_{p,r}},
\end{align*}
which gives the desired result of Proposition \ref{pro5-1}.
\end{proof}
The following proposition plays a key role in the proof of Theorem \ref{th4}.
\begin{proposition}\label{pro5-2}
Under the assumptions of Theorem \ref{th4}, we have
\bbal
\|\mathbf{S}_{t}(S_nu_0+v^{n,\omega}_{0,m_n})-\mathbf{S}_{t}(S_nu_0)
-\mathbf{S}_{t}(v^{n,\omega}_{0,m_n})\|_{B^{s}_{p,r}}\leq C2^{n(1-\theta)}e^{2^n\theta}\mathbf{F}_n^{\theta}\quad\text{with}\quad \theta=\fr{1}{s+1},
\end{align*}
where we denote
$$\mathbf{F}_n=\int_0^t\big\|\mathbf{S}_{t}(S_nu_0)\mathbf{S}_{t}(v^{n,\omega}_{0,m_n}),\pa^2_{xx}(\mathbf{S}_{t}(S_nu_0)\mathbf{S}_{t}(v^{n,\omega}_{0,m_n}))
,\pa_x\mathbf{S}_{t}(S_nu_0)\pa_x\mathbf{S}_{t}(v^{n,\omega}_{0,m_n})\big\|_{L^p}\dd \tau.$$
\end{proposition}
\begin{proof} For the sake of simplicity, we set $\mathbf{w}=\mathbf{S}_{t}(S_nu_0+v^{n,\omega}_{0,m_n})-\mathbf{S}_{t}(S_nu_0)-\mathbf{S}_{t}(v^{n,\omega}_{0,m_n})$.
By the interpolation inequality (see Lemma \ref{le2}), we obtain
\bal\label{lyz}
\|\mathbf{w}\|_{B^s_{p,r}}\leq C\|\mathbf{w}\|^{\theta}_{B^0_{p,\infty}}\|\mathbf{w}\|^{1-\theta}_{B^{s+1}_{p,\infty}}
\leq C2^{n(1-\theta)}\|\mathbf{w}\|^{\theta}_{L^p}.
\end{align}
Next, we need to estimate the $L^p$-norm of $\mathbf{w}$. It is obvious that $\mathbf{w}$ solves
\begin{align}\label{m0}
\begin{cases}
\partial_t\mathbf{w}+\mathbf{S}_{t}(S_nu_0+v^{n,\omega}_{0,m_n})\pa_x\mathbf{w}=
-\mathbf{w}\pa_x\left(\mathbf{S}_{t}(S_nu_0)+\mathbf{S}_{t}(v^{n,\omega}_{0,m_n})\right)-\sum_{i=1}^3F_{i}
,\\
\mathbf{w}(0,x)=0,
\end{cases}
\end{align}
where
\bbal
&F_{1}=\pa_x\big(\mathbf{S}_{t}(S_nu_0)\mathbf{S}_{t}(v^{n,\omega}_{0,m_n})\big),\\
&F_{2}=2\mathcal{B}(\mathbf{S}_{t}(S_nu_0),\mathbf{S}_{t}(v^{n,\omega}_{0,m_n})),\\
&F_{3}=\mathcal{B}(\mathbf{w}, \mathbf{S}_{t}(S_nu_0+v^{n,\omega}_{0,m_n})+\mathbf{S}_{t}(S_nu_0)+\mathbf{S}_{t}(v^{n,\omega}_{0,m_n})).
\end{align*}
Taking the inner product of Eq. $\eqref{m0}_1$ with $|\mathbf{w}|^{p-2}\mathbf{w}$ with $p\in[1,\infty)$, we obtain
\begin{align}\label{z1}
\fr1p\frac{\dd}{\dd t}\|\mathbf{w}\|^p_{L^p}&=-\int_{\R}\pa_x\left(\mathbf{S}_{t}(S_nu_0)+\mathbf{S}_{t}(v^{n,\omega}_{0,m_n})\right)|\mathbf{w}|^{p}\dd x\nonumber\\
&\quad+\fr1p\int_{\R}\pa_x\mathbf{S}_{t}(S_nu_0+v^{n,\omega}_{0,m_n})|\mathbf{w}|^{p}\dd x-\int_{\R}\sum_{i=1}^3F_{i}|\mathbf{w}|^{p-2}\mathbf{w}\dd x\nonumber\\
&\les \left(\|\pa_x(\mathbf{S}_{t}(S_nu_0)+\mathbf{S}_{t}(v^{n,\omega}_{0,m_n}))\|_{L^\infty}
+\|\pa_x\mathbf{S}_{t}(S_nu_0+v^{n,\omega}_{0,m_n})\|_{L^\infty}\right)\|\mathbf{w}\|^p_{L^p}\nonumber\\
&\quad
+\big\|F_{1},F_{2},F_{3}\big\|_{L^p}\|\mathbf{w}\|^{p-1}_{L^p}\nonumber\\
&\les \|\mathbf{w}\|^p_{L^p}
+\|\pa_x\mathbf{w}\|_{L^p}\|\mathbf{w}\|^{p-1}_{L^p}
+\|F_{1},F_{2}\|_{L^p}\|\mathbf{w}\|^{p-1}_{L^p},
\end{align}
where we have used the fact
$\|\mathcal{B}(f,g)\|_{L^p}\les\|fg\|_{L^p}+\|\pa_xf\pa_xg\|_{L^p}$.

Then \eqref{z1} reduces to
\begin{align}\label{ch1}
\frac{\dd}{\dd t}\|\mathbf{w}\|_{L^p}
&\les\|\mathbf{w},\pa_x\mathbf{w}\|_{L^p}+\|F_{1},F_{2}\|_{L^p}.
\end{align}

Next, to close \eqref{ch1}, we have to estimate the $L^p$-norm of $\pa_x\mathbf{w}$.

Setting $\mathbf{v}=\pa_x\mathbf{w}$ with $\mathbf{v}(0,x)=0$, then we have
\begin{align}\label{m}
\partial_t\mathbf{v}+\mathbf{S}_{t}(S_nu_0+v^{n,\omega}_{0,m_n})\pa_x\mathbf{v}&=
-\mathbf{v}\pa_x\left(\mathbf{S}_{t}(S_nu_0+v^{n,\omega}_{0,m_n})+\mathbf{S}_{t}(S_nu_0)+\mathbf{S}_{t}(v^{n,\omega}_{0,m_n})\right)\nonumber\\
&\quad-\mathbf{w}\pa_{xx}^2\left(\mathbf{S}_{t}(S_nu_0)+\mathbf{S}_{t}(v^{n,\omega}_{0,m_n})\right)-\sum_{i=1}^3\pa_xF_{i}.
\end{align}
Taking the inner product of Eq. \eqref{m} with $|\mathbf{v}|^{p-2}\mathbf{v}$, we obtain
\begin{align}\label{z2}
\fr1p\frac{\dd}{\dd t}\|\mathbf{v}\|^p_{L^p}
&=-\int_{\R}\sum_{i=1}^3\pa_xF_{i}|\mathbf{v}|^{p-2}\mathbf{v}\dd x-\int_{\R}\pa_x\left(\mathbf{S}_{t}(S_nu_0+v^{n,\omega}_{0,m_n})+\mathbf{S}_{t}(S_nu_0)+\mathbf{S}_{t}(v^{n,\omega}_{0,m_n})\right)
|\mathbf{v}|^{p}\dd x\nonumber\\
&\quad+\fr1p\int_{\R}\pa_x\mathbf{S}_{t}(S_nu_0+v^{n,\omega}_{0,m_n})|\mathbf{v}|^{p}\dd x
-\int_{\R}\mathbf{w}\pa_{xx}^2\left(\mathbf{S}_{t}(S_nu_0)+\mathbf{S}_{t}(v^{n,\omega}_{0,m_n})\right)|\mathbf{v}|^{p-2}\mathbf{v}\dd x\nonumber\\
&\les\big\|\sum_{i=1}^3\pa_xF_{i}\big\|_{L^p}\|\mathbf{v}\|^{p-1}_{L^p}+\|\pa_x\mathbf{S}_{t}(S_nu_0),\pa_x\mathbf{S}_{t}(v^{n,\omega}_{0,m_n}),\pa_x\mathbf{S}_{t}(S_nu_0+v^{n,\omega}_{0,m_n})\|_{L^\infty}
\|\mathbf{v}\|^p_{L^p}\nonumber\\
&\quad+\left\|\pa_{xx}^2\left(\mathbf{S}_{t}(S_nu_0)+\mathbf{S}_{t}(v^{n,\omega}_{0,m_n})\right)\right\|_{L^\infty}\|\mathbf{w}\|_{L^p}\|\mathbf{v}\|^{p-1}_{L^p}
\nonumber\\
&\les\|\mathbf{v}\|^p_{L^p}+2^n\|\mathbf{w}\|_{L^p}\|\mathbf{v}\|^{p-1}_{L^p}
+\big\|\pa_xF_{1},\pa_xF_{2}\big\|_{L^p}\|\mathbf{v}\|^{p-1}_{L^p},
\end{align}
where we have used the fact due to $s-1>\max\{\frac12,\frac1p\}$
$$\left\|\pa_{xx}^2\left(\mathbf{S}_{t}(S_nu_0)+\mathbf{S}_{t}(v^{n,\omega}_{0,m_n})\right)\right\|_{L^\infty}\les \left\|\pa_{xx}^2\left(\mathbf{S}_{t}(S_nu_0)+\mathbf{S}_{t}(v^{n,\omega}_{0,m_n})\right)\right\|_{B^{s-1}_{p,r}}\les2^n.$$
Then \eqref{z2} reduces to
\begin{align}\label{ch2}
\frac{\dd}{\dd t}\|\mathbf{v}\|_{L^p}
&\les2^n\|\mathbf{w}\|_{L^p}+\|\mathbf{v}\|_{L^p}+\big\|\pa_xF_{1},\pa_xF_{2}\big\|_{L^p}.
\end{align}

Combining \eqref{ch1} and \eqref{ch2} yields that
\begin{align}\label{ch3}
\|\mathbf{w},\pa_x\mathbf{w}\|_{L^p}&\les e^{2^nt}\int_0^t\sum_{i=1}^2\big\|F_{i},\pa_xF_{i}\big\|_{L^p}\dd\tau
\les e^{2^n}\mathbf{F}_n.
\end{align}
Particularly, we would like to emphasize that \eqref{ch3} holds for $p=\infty$ (just letting $p\rightarrow\infty$ in \eqref{ch3}).

Inserting \eqref{ch3} into \eqref{lyz}, then we complete the proof of Proposition \ref{pro5-2}.
\end{proof}

With Propositions \ref{pro1}--\ref{pro2} in hand, we can  prove Theorem \ref{th4} by dividing it into three steps.

{\bf Step 1.}\; By Theorem \ref{th0}, we can construct two sequences of initial data $\mathbf{f}_n$ and $\mathbf{f}_n+\mathbf{g}_n$ (see \eqref{fg}) which satisfy
$\|\mathbf{f}_n\|_{B^s_{p,r}}\lesssim 1$ and $\lim_{n\rightarrow \infty}\|\mathbf{g}_n\|_{B^s_{p,r}}= 0
$
and find two sequences of solutions $\mathbf{S}_{t}(\mathbf{f}_n+\mathbf{g}_n)$ and $\mathbf{S}_{t}(\mathbf{f}_n)$ which satisfy $\|\mathbf{S}_{t}(\mathbf{f}_n+\mathbf{g}_n),\mathbf{S}_{t}(\mathbf{f}_n)\|_{B^s_{p,r}}\leq C$,
but for a short time $t\in[0,T_0]$
\bbal
\liminf_{n\rightarrow \infty}\|\mathbf{S}_t(\mathbf{f}_n+\mathbf{g}_n)-\mathbf{S}_t(\mathbf{f}_n)\|_{B^s_{p,r}}\geq c_0t>0.
\end{align*}

{\bf Step 2.}\;
Notice that
$\mathbf{S}_{t}(v^{n,\omega}_{0,m_n})=\mathbf{S}_{t}(\mathbf{f}_n+\omega \mathbf{g}_n)(t,x-m_n):=u^{\omega}_n(t,x-m_n)$, for fixed $n$ and any $(t,x)$, thanks to the smoothness and decay, we have
\bbal
&\lim_{|y|\rightarrow \infty} \mathbf{S}_{t}(S_nu_0(x))\mathbf{S}_{t}(\mathbf{f}_n+\omega \mathbf{g}_n)(t,x-y)=0,\\
&\lim_{|y|\rightarrow \infty} \pa^2_{xx}(\mathbf{S}_{t}(S_nu_0(x))\mathbf{S}_{t}(\mathbf{f}_n+\omega \mathbf{g}_n)(t,x-y))=0,\\
&\lim_{|y|\rightarrow \infty} \pa_x\mathbf{S}_{t}(S_nu_0(x))\pa_x\mathbf{S}_{t}(\mathbf{f}_n+\omega \mathbf{g}_n)(t,x-y)=0.
\end{align*}
Also, we have
\bbal
&|\mathbf{S}_{t}(S_nu_0(x))\mathbf{S}_{t}(\mathbf{f}_n+\omega \mathbf{g}_n)(t,x-y)|\leq M|\mathbf{S}_{t}(S_nu_0(x))|,\\
&|\mathbf{S}_{t}(S_nu_0(x))\pa^2_{xx}\mathbf{S}_{t}(\mathbf{f}_n+\omega \mathbf{g}_n)(t,x-y)|\leq M|\mathbf{S}_{t}(S_nu_0(x))|,\\
&|\pa_x\mathbf{S}_{t}(S_nu_0(x))\pa_x\mathbf{S}_{t}(\mathbf{f}_n+\omega \mathbf{g}_n)(t,x-y)|\leq M|\pa_x\mathbf{S}_{t}(S_nu_0(x))|,\\
&|\pa^2_{xx}\mathbf{S}_{t}(S_nu_0(x))\mathbf{S}_{t}(\mathbf{f}_n+\omega \mathbf{g}_n)(t,x-y)|\leq M|\pa^2_{xx}(\mathbf{S}_{t}(S_nu_0(x)))|.
\end{align*}
By the Lebesgue Dominated Convergence Theorem, one has for $\omega=0,1$
\bbal
&\lim_{|y|\rightarrow \infty} \int_0^t\big\|\mathbf{S}_{t}(S_nu_0)u^{\omega}_n(t,x-y),\pa^2_{xx}(\mathbf{S}_{t}(S_nu_0)u^{\omega}_n(t,x-y))
,\pa_x\mathbf{S}_{t}(S_nu_0)\pa_xu^{\omega}_n(t,x-y)\big\|_{L^p}\dd \tau=0.
\end{align*}
By Proposition \ref{pro5-2}, for fixed $n$ and any $(t,x)$, one can find $m_n$ with $|m_n|$ sufficiently large such that
\bbal
\|\mathbf{S}_{t}(S_nu_0+v^{n,\omega}_{0,m_n})-\mathbf{S}_{t}(S_nu_0)
-\mathbf{S}_{t}(v^{n,\omega}_{0,m_n})\|_{B^{s}_{p,r}}\leq \fr{c_0t}{4}, \quad \omega=0,1.
\end{align*}

{\bf Step 3.}\; We decompose the difference of $\mathbf{S}_{t}(u_0+v^{n,1}_{0,m_n})$ and $\mathbf{S}_{t}(u_0+v^{n,0}_{0,m_n})$ as follows
\bbal
\mathbf{S}_{t}(&u_0+v^{n,1}_{0,m_n})-\mathbf{S}_{t}(u_0+v^{n,0}_{0,m_n})=\mathbf{S}_{t}(v^{n,1}_{0,m_n})-\mathbf{S}_{t}(v^{n,0}_{0,m_n})\\
&\quad+\left(\underbrace{\mathbf{S}_{t}(u_0+v^{n,1}_{0,m_n})-\mathbf{S}_{t}(S_nu_0+v^{n,1}_{0,m_n})}_{:=\mathbf{I}_1}+\underbrace{\mathbf{S}_{t}(S_nu_0+v^{n,1}_{0,m_n})-\mathbf{S}_{t}(S_nu_0)
-\mathbf{S}_{t}(v^{n,1}_{0,m_n})}_{:=\mathbf{I}_2}\right)\\
&\quad-\left(\underbrace{\mathbf{S}_{t}(u_0+v^{n,0}_{0,m_n})-\mathbf{S}_{t}(S_nu_0+v^{n,0}_{0,m_n})}_{:=\mathbf{I}_3}+\underbrace{\mathbf{S}_{t}(S_nu_0+v^{n,0}_{0,m_n})-\mathbf{S}_{t}(S_nu_0)
-\mathbf{S}_{t}(v^{n,0}_{0,m_n})}_{:=\mathbf{I}_4}\right).
\end{align*}
Hence, from Proposition \ref{pro5-1} and {\bf Step 2} we deduce
\bbal
\|\mathbf{S}_{t}(u_0&+v^{n,1}_{0,m_n})-\mathbf{S}_{t}(u_0+v^{n,0}_{0,m_n})\|_{B^s_{p,r}}
\geq \|\mathbf{S}_{t}(v^{n,1}_{0,m_n})-\mathbf{S}_{t}(v^{n,0}_{0,m_n})\|_{B^s_{p,r}}-\sum_{i=1}^4\|\mathbf{I}_i\|_{B^s_{p,r}}\\
&\geq \|\mathbf{S}_{t}(v^{n,1}_{0,m_n})-\mathbf{S}_{t}(v^{n,0}_{0,m_n})\|_{B^s_{p,r}}-C\|(\mathrm{Id}-S_n)u_0\|_{B^{s}_{p,r}}-\fr{c_0t}{2},
\end{align*}
which follows from {\bf Step 1} that for $t$ small enough
\bbal
\liminf_{n\rightarrow \infty}\|\mathbf{S}_{t}(u_0+v^{n,1}_{0,m_n})-\mathbf{S}_{t}(u_0+v^{n,0}_{0,m_n})\|_{B^s_{p,r}}\geq \fr{c_0t}{2}.
\end{align*}
We also notice that
\bbal
\lim_{n\to\infty}\|v^{n,1}_{0,m_n}-v^{n,0}_{0,m_n}\|_{B^s_{p,r}}=\lim_{n\to\infty}\|g_{n}(\cdot-m_n)\|_{B^s_{p,r}}=0.
\end{align*}
This completes the proof of Theorem \ref{th4}. {\hfill $\square$}

\section*{Acknowledgments}
J. Li is supported by the National Natural Science Foundation of China (12161004), Training Program for Academic and Technical Leaders of Major Disciplines in Ganpo Juncai Support Program(20232BCJ23009) and Jiangxi Provincial Natural Science Foundation (20224BAB201008). Y. Yu is supported by the National Natural Science Foundation of China (12101011). W. Zhu is supported by the National Natural Science Foundation of China (12201118) and Guangdong Basic and Applied Basic Research Foundation (2021A1515111018).

\section*{Declarations}
\noindent\textbf{Data Availability} No data was used for the research described in the article.

\vspace*{1em}
\noindent\textbf{Conflict of interest}
The authors declare that they have no conflict of interest.

\end{document}